\documentclass[10pt]{amsart}

\usepackage{amsmath,amsthm,amsfonts,amssymb}
\usepackage{color}
\usepackage{verbatim}
\usepackage{pictexwd}
\usepackage{graphicx}
\usepackage{hyperref}
\usepackage[all,dvips,arc,curve,color,frame]{xy}
\usepackage{algorithm}
\usepackage{algorithmicx}
\usepackage{algpseudocode}

\theoremstyle{definition}

\newtheorem*{theorem*}{Theorem}
\newtheorem*{proposition*}{Proposition}

\newtheorem{alg}{Algorithm}

\algnewcommand\algorithmicinitialsetup{\textbf{Initial Setup:}}
\algnewcommand\InitialSetup{\item[\algorithmicinitialsetup]}
\algnewcommand\algorithmicalicepublickey{\textbf{Alice's Public Key:}}
\algnewcommand\AlicePublicKey{\item[\algorithmicalicepublickey]}
\algnewcommand\algorithmicbobpublickey{\textbf{Bob's Public Key:}}
\algnewcommand\BobPublicKey{\item[\algorithmicbobpublickey]}
\algnewcommand\algorithmicaliceprivatekey{\textbf{Alice's Private Key:}}
\algnewcommand\AlicePrivateKey{\item[\algorithmicaliceprivatekey]}
\algnewcommand\algorithmicbobprivatekey{\textbf{Bob's Private Key:}}
\algnewcommand\BobPrivateKey{\item[\algorithmicbobprivatekey]}
\algnewcommand\algorithmicsharedkey{\textbf{Shared Key:}}
\algnewcommand\SharedKey{\item[\algorithmicsharedkey]}
\algnewcommand\algorithmicalicesharedkey{\textbf{Alice's Shared Key:}}
\algnewcommand\AliceSharedKey{\item[\algorithmicalicesharedkey]}
\algnewcommand\algorithmicbobsharedkey{\textbf{Bob's Shared Key:}}
\algnewcommand\BobSharedKey{\item[\algorithmicbobsharedkey]}
\algnewcommand\algorithmicencoding{\textbf{Encryption:}}
\algnewcommand\Encoding{\item[\algorithmicencoding]}
\algnewcommand\algorithmicdecoding{\textbf{Decryption:}}
\algnewcommand\Decoding{\item[\algorithmicdecoding]}

\def\ova{{\overline{a}}}
\def\ovb{{\overline{b}}}

\def\ove{{\overline{e}}}
\def\ovo{{\overline{0}}}
\def\ovs{{\overline{s}}}
\def\ovt{{\overline{t}}}

\def\ovv{{\overline{v}}}

\def\CO{{\mathcal O}}

\newcommand{\gp}[1]{{\left\langle #1 \right\rangle}}

\newcommand{\rb}[1]{{\left( #1 \right)}}

\def\MQ{{\mathbb{Q}}}
\def\MZ{{\mathbb{Z}}}

\let\oldmarginpar\marginpar
\renewcommand\marginpar[1]{\-\oldmarginpar[\raggedleft\footnotesize #1]%
{\raggedright\footnotesize #1}}

\title{Analysis of a certain polycyclic-group-based cryptosystem}
\date{\today}

\author[M. Kotov]{Matvei Kotov}
\author[A. Ushakov]{Alexander Ushakov}
\address{Department of Mathematics, Stevens Institute of Technology, Hoboken, NJ, USA}
\email{mkotov,aushakov@stevens.edu}
\thanks{The second author has been partially supported by NSA Mathematical
Sciences Program grant number H98230-14-1-0128}

\begin{document}

\begin{abstract}
We investigate security properties
of the Anshel-Anshel-Goldfeld commutator key-establishment protocol \cite{AAG} used with
certain polycyclic groups described in \cite{Eick}.
We show that despite low success of the length based attack shown in
\cite{Garber-Kahrobaei-Lam:2015} the protocol can be broken by a deterministic
polynomial-time algorithm.

\end{abstract}

\maketitle

\section{Introduction}

In this paper we analyze  the \textbf{commutator key-establishment} protocol \cite{AAG}
used with certain polycyclic groups described in \cite{Eick}.
The commutator key-establishment (CKE) protocol is a two-party protocol
performed as follows.
\begin{itemize}
\item
Fix a group $G$ (called \textbf{the platform group}) and a set of generators $g_1,\ldots,g_k$
for $G$. All this information is made public.
\item
Alice prepares a tuple of elements $\ova = (a_1,\ldots,a_{N_1})$ called
\textbf{Alice's public tuple}.
Each $a_i$ is generated randomly as a product of $g_i$'s and their inverses.
\item
Bob prepares a tuple of elements $\ovb=(b_1,\ldots,b_{N_2})$ called the
\textbf{Bob's public tuple}.
Each $b_i$ is generated randomly as a product of $g_i$'s and their inverses.
\item
Alice generates a random element $A$ as a product
$a_{s_1}^{\varepsilon_1} \ldots a_{s_L}^{\varepsilon_L}$
of $a_i$'s and their inverses. The element $A$ (or more precisely its factorization)
is called the \textbf{Alice's private element}.
\item
Bob generates a random element $B$ as a product
$b_{t_1}^{\delta_1} \ldots a_{t_L}^{\delta_L}$
of $b_i$'s and their inverses, called the \textbf{Bob's private element}.
\item
Alice publishes the tuple of conjugates $\ovb^A = (A^{-1}b_1A,\ldots,A^{-1}b_{N_2}A)$.
\item
Bob publishes the tuple of conjugates $\ova^B = (B^{-1}a_1B,\ldots,B^{-1}a_{N_1}B)$.
\item
Finally, Alice computes the element $K_A$ as a product:
$$
A^{-1} \cdot \rb{B^{-1}a_{s_1}^{\varepsilon_1}B \ldots B^{-1}a_{s_L}^{\varepsilon_L}B}
$$
using the elements of Bob's conjugate tuple $\ova^B$.
\item
Bob computes the key $K_B$ as a product:
$$
\rb{A^{-1}b_{t_1}^{\delta_1}A \ldots A^{-1}b_{t_L}^{\delta_L}A}^{-1} \cdot B
$$
using the elements of Alice's conjugate tuple $\ovb^A$.
\end{itemize}
It is easy to check that $K_A=K_B = A^{-1}B^{-1}AB$ in $G$.
The obtained commutator is the \textbf{shared key}.

Security of the commutator key establishment protocol is based on
computational hardness of computing the commutator $[A,B]$
based on the intercepted public information -- the tuples $\ova,\ovb$
and their conjugates $\ova^B,\ovb^A$. In practice it is often achieved by
solving systems of conjugacy equations for $A$ and $B$, i.e., finding
$X=A'$ and $Y=B'$ satisfying:
$$
\left\{
\begin{array}{l}
X^{-1} b_1 X = b_1',  \\
\ldots \\
X^{-1} b_{N_1} X = b_{N_1}',
\end{array}
\right.
\quad\mbox{ and }\quad
\left\{
\begin{array}{l}
Y^{-1} a_1 Y = a_1',  \\
\ldots \\
Y^{-1} a_{N_2} Y = a_{N_1}',
\end{array}
\right.
$$
and computing $K'=[A',B']$. In general it can happen that $K'\ne K$
as explained in \cite{SU2}, but as practice shows very often $K=K'$
(for instance, as in \cite{Hofheinz-Steinwandt:2003}).

A big advantage of the commutator key-establishment protocol
over other group-based protocols
is that it can be used with any group $G$ satisfying certain computational
properties.
Originally, the group of braids $B_n$ was suggested to use as a platform group,
but after a series of attacks it became clear that $B_n$ can not provide good security.
But the search for a good group is still very active and in \cite{Eick}
a certain class of polycyclic groups was proposed to be used with CKE.
In this paper we show that that class can not provide good security.
For more on group-based cryptography see \cite{MSU_book:2011}.

\subsection{Outline}
In Section \ref{se:platform} we define the class of groups
under investigation and discuss two different ways to represent the elements.
In Sections \ref{se:semidirect_attack} and
\ref{se:polycyclic_attack} we describe the attacks on different group presentations.

\section{The platform group}
\label{se:platform}

Consider an irreducible monic polynomial $f(x)\in\MZ[x]$ and define a field:
$$
F = \MQ[x]/(f).
$$
The \textbf{ring of integers} of $F$ is defined as:
$$
\CO_F = \{a\in F \mid a \mbox{ is a zero of a monic polynomial } g(x)\in\MZ[x] \}
$$
and its \textbf{group of units}:
$$
U_F = \{a \mid a^{-1} \in \CO_F\}.
$$
A semidirect product $U_F \ltimes \CO_F$ of $U_F$ and $\CO_F$ is defined as a
Cartesian product $U_F \times \CO_F$ equipped with the following binary operation:
\begin{equation}\label{eq:semi_prod}
(\alpha,a)\cdot(\beta,b) = (\alpha\beta, a\beta + b).
\end{equation}
The constructed group $G_F$ is the platform group in \cite{Eick}.
It is easy to see that $G_F$ is polycyclic and metabelian and
there are several different ways to represent $G_F$.
\begin{itemize}
\item[(a)]
One can work with $G_F$ as it is defined above, i.e., as a semidirect product,
in which case its elements are represented as pairs and
multiplication (\ref{eq:semi_prod}) is used.
\item[(b)]
One can construct a polycyclic presentation for $G_F$ and work with
its elements as with words over the generating set.
\end{itemize}
Unfortunately, neither \cite{Eick} nor \cite{Garber-Kahrobaei-Lam:2015}
give any detail on how to treat $G_F$.
Since computational properties of the same group can vary depending on a way
we represent its elements,
in the next sections we discuss both presentations of $G_F$.

\subsection{$G_F$ as a set of pairs of matrices}
\label{se:semidirect}

There are different ways to represent the elements of $F$.
For instance, elements in $F$ can be represented as polynomials
over $\MQ$ of degree up to $n-1$ with addition and multiplication performed
modulo the original polynomial $f$. Also one can represent elements in $F$
by matrices as described below.
Recall that the \textbf{companion matrix} for a monic polynomial
$f = x^n+c_{n-1}x^{n-1}+\ldots+c_1x+c_0$
is a matrix of the form:
$$
M=
\left[
\begin{array}{ccccc}
0 & 0 & \ldots & 0 & -c_0 \\
1 & 0 & \ldots & 0 & -c_1 \\
0 & 1 & \ldots & 0 & -c_2 \\
\vdots & \vdots & \ddots & \vdots & \vdots \\
0 & 0 & \ldots & 1 & -c_{n-1}
\end{array}
\right]
$$
The characteristic  and minimal polynomial of $M$ is $f$ and the set of matrices:
\begin{equation}\label{eq:field_F}
F = \left\{a_0E + a_1M + a_2M^2 + \ldots + a_{n-1}M^{n-1} \,|\, a_0, \ldots, a_{n-1} \in \mathbb{Q}\right\}.
\end{equation}
equipped with the usual matrix addition and multiplication is the field $F$.
The correspondence between two presentations is obvious:
$$
a_0+a_1x+\ldots +a_{n-1}x^{n-1} \longleftrightarrow
a_0+a_1M+\ldots +a_{n-1}M^{n-1}
$$
and choosing a particular presentation we do not change computational properties
of $F$. Here we choose matrix presentation for $F$.

Let $O_1, \ldots, O_n$ be a basis of the ring of integers $\CO_F$,
where each $O_i$ is a matrix. Hence:
$$\CO_F = \left\{a_1O_1 + a_2O_2 + \ldots + a_nO_n \,|\, a_1, \ldots, a_n \in \mathbb{Z}\right\}.$$
Let $\{U_1, \ldots, U_m\}$ be a generating set for the group $U_F$,
where every $U_i$ is a matrix. Hence:
$$
U_F = \left\{U_1^{a_1}\cdot U_2^{a_2}\cdot\ldots\cdot U_m^{a_m} \,|\, a_1, \ldots, a_m \in \mathbb{Z}\right\}.
$$
By Dirichlet theorem \cite[Chapter 8]{holt2005handbook}
$U_F \cong \mathbb{Z}_k\times \mathbb{Z}^{m-1}$, where $m = s + t - 1$, $s$ is the number of real field monomorphisms $F \to \mathbb{R}$, and $2t$ is the number of complex field monomorphisms $F \to \mathbb{C}$. Without loss of generality it can be assumed that $U_1^k = E$.

Now naturally the group $G_F = U_F \ltimes \CO_F$ is a set of pairs of matrices:
$$G = \{(C, S) \,|\, C \in U_F, S \in \CO_F\},$$
equipped with multiplication given by:
\begin{equation}\label{eq:mul}
(C,S)\cdot(D,T) = (CD, SD + T).
\end{equation}
It is easy to check that the inverse in $U_F \ltimes \CO_F$ can be computed as
\begin{equation}\label{eq:inv}
(C,S)^{-1} = (C^{-1}, -SC^{-1}),
\end{equation}
which gives the following expression for the conjugate of $(B,T)$ by $(C,S)$
\begin{equation}\label{eq:conj}
(D, T)^{(C, S)} = (C,S)^{-1}(D,T)(C,S) = (D, S(E-D) + TC),
\end{equation}
where $E$ is the identity matrix.

\subsection{$G_F$ given by polycyclic presentation}
\label{se:poly_pres}

Recall that a group $G$ is called polycyclic if there exists a subnormal series of $G$:
$$
G = G_0 \rhd G_2 \rhd G_3 \ldots \rhd G_n = \{1\},
$$
with cyclic factors $G_{i-1}/G_{i}$.
Denote $[G_{i-1}:G_{i}]$ by $r_i$ and put $I=\{i \mid r_i<\infty\}$.
Relative to the series above one can
find a generating set $g_1,\ldots,g_n$ for $G$ satisfying $\gp{G_{i},g_i} = G_{i-1}$.
Every element $g\in G$ can be uniquely expressed as a product
$g = g_1^{e_1} \ldots g_n^{e_n},$ where $e_i \in \MZ$, $i=1,\ldots, n$, and
$0\le e_i<r_i$ if  $i\in I$.
The polycyclic group $G$ has a finite presentation of the form:
\begin{equation}\label{eq:wordpres}
G =
 \left\langle
\begin{array}{lcl}g_1,\ldots,g_n & \bigg{|} &
\begin{array}{ll}
 g_j^{g_i} = w_{ij},
g_j^{g_i^{-1}} = v_{ij} & \textrm{for } 1\le i < j\le n, \\
 g_k^{r_k}=u_{k} & \textrm{for } k\in I
\end{array}
\end{array}
\right\rangle,
\end{equation}
where $w_{ij}$, $v_{ij}$, and $u_i$ are words in $g_{i+1},\ldots,g_n$. This
presentation is called a \textbf{polycyclic presentation}.
For more details see \cite[Chapter 8]{holt2005handbook}.

It is straightforward to find a polycyclic presentation for
the group $G_F = U_F \ltimes \CO_F$.
It has generators $g_1, \ldots, g_m, g_{m+1}, \ldots, g_{m+n}$, where
$g_1, \ldots, g_m$ correspond to the pairs
$(U_1, O), \ldots, (U_m, O) \in U_F \ltimes \CO_F$
($O$ is the zero matrix),
and $g_{m+1}, \ldots, g_{m+n}$ correspond to the pairs
$(E, O_1), \ldots, (E, O_n) \in U_F \ltimes \CO_F$
($E$ is the identity matrix).
The set of relations for $G$ is formed as follows.
\begin{itemize}
\item
$g_{m+j}^{g_i} = g_{m+1}^{a_{ij1}}\ldots g_{m+n}^{a_{ijn}}$,
$i = 1,\ldots, m$, $j=1, \ldots, n$, and $a_{ij1}, \ldots, a_{ijn}$
are the coefficients in the expression $O_jU_i = a_{ij1}O_1 + \ldots + a_{ijn}O_n$,
\item
$g_{m+j}^{g_i^{-1}} = g_{m+1}^{b_{ij1}}\ldots g_{m+n}^{b_{ijn}}$,
$i = 1,\ldots, m$, $j=1, \ldots, n$, and $a_{ij1}, \ldots, b_{ijn}$
are the coefficients in the expression $O_jU_i^{-1} = b_{ij1}O_1 + \ldots + b_{ijn}O_n$,
\item $g_1^k = e$,
\item $[g_i, g_j] = e$, $1 \leq i < j \leq m$,
\item $[g_i, g_j] = e$, $m+1 \leq i < j \leq m+n$.
\end{itemize}

\section{Attack on semidirect product}
\label{se:semidirect_attack}

In this section we assume that the group $G_F$ is given as a semidirect product
and the field $F$ is described using matrices as in (\ref{eq:field_F}).
The general idea behind the attack is to extend
the group $G_F$ and work in $G_F^\ast=F^\ast \ltimes F$.
The group $G_F^\ast$
is, in general, not finitely generated and hence is not polycyclic.
Nevertheless the elements of $G^\ast$ can be effectively represented by pairs of matrices
as described in Section \ref{se:semidirect}.

Consider a system of conjugacy equations related to the Alice's private key:
\begin{equation}\label{eq:conds}
\left\{
\begin{array}{ll}
X^{-1}b_1X &=b_1', \\
&\ldots \\
X^{-1}b_{N_2}X &=b'_{N_2},\\
\end{array}
\right.
\end{equation}
with unknown $X\in U_F\ltimes \CO_F$.
We treat the system as a system over $F^\ast \ltimes F$ and hence:
$$
X=(C,S),\
b_i=(B_i, T_i),\
b_i'=(B_i', T_i')\
\mbox{ in } F^\ast \ltimes F.
$$
Using (\ref{eq:conj}) we get the following system
of $N_2$ linear equations over the field $F$ with two unknowns $C$ and $S$:
\begin{equation}\label{eq:thesys}
\left\{
\begin{array}{lcl}
S(E - B_1) + T_1C & = & T'_1, \\
& \vdots & \\
S(E - B_{N_2}) + T_{N_2}C & = &T'_{N_2}.
\end{array}
\right.
\end{equation}
It has a unique solution when the coefficient matrix of the system has rank $2$ over the field $F$, in which case the obtained solution $A'$ is the same as the original Alice's private key.
We call the described approach ``field based attack'' or simply \textbf{FBA}.

The described attack was implemented in GAP \cite{GAP4}.
Its implementation can be found in \cite{github:Matvej}.
The table below compares success rate and time efficiency
of our attack and the attack in~\cite{Garber-Kahrobaei-Lam:2015}.
Our tests were run on Intel Core i5 1.80GHz computer with 4GB of RAM, Ububtu 12.04, GAP 4.7.

\begin{center}
\begin{tabular}{|c|c|c|c|c|c|c|c|}
\hline
\raisebox{-4ex}[0cm][0cm]{Polynomial} & \raisebox{-4ex}[0cm][0cm]{$h(G)$} & \multicolumn{2}{|c|}{LBA w/ dynamic} & \multicolumn{2}{|c|}{\raisebox{-1.5ex}[0cm][0cm]{FBA, $L=5$}}  & \multicolumn{2}{|c|}{\raisebox{-1.5ex}[0cm][0cm]{FBA, $L=100$}} \\
 & & \multicolumn{2}{|c|}{set, $L=5$} & \multicolumn{2}{|c|}{}  & \multicolumn{2}{|c|}{} \\
\cline{3-8}
& & \raisebox{-1.5ex}[0cm][0cm]{Time} & Success & \raisebox{-1.5ex}[0cm][0cm]{Time}   &  Success & \raisebox{-1.5ex}[0cm][0cm]{Time} & Success\\
& & & rate  &    &  rate & & rate\\
\hline
 $x^2-x-1$      & 3  &  0.20 h & 100\% &   2.4 s & 100\% &   2.8 s & 100\% \\
 $x^5-x^3-1$    & 7  & 76.87 h &  35\% &   3.4 s & 100\% &   5.3 s & 100\% \\
 $x^7-x^3-1$    & 10 & 94.43 h &   8\% &   5.2 s & 100\% &   9.7 s & 100\% \\
 $x^9-7x^3-1$   & 14 & 95.18 h &   5\% &  23.1 s & 100\% &  57.7 s & 100\% \\
 $x^{11}-x^3-1$ & 16 & 95.05 h &   5\% &  15.3 s & 100\% &  29.5 s & 100\% \\
 $x^{15}-x-2$   & 22 &      -- &    -- & 694.8 s & 100\% & 607.4 s & 100\% \\
 $x^{20}-x-1$   & 30 &      -- &    -- & 208.5 s & 100\% & 192.8 s & 100\% \\
\hline
\end{tabular}
\end{center}

The first four columns of this table are taken from~\cite{Garber-Kahrobaei-Lam:2015}.
For our tests we used the same parameter values: $N_1 = N_2 = 20$, and the same
number of tests: $100$.

\section{Attack on polycyclic presentation}
\label{se:polycyclic_attack}

In this section we assume that $G_F$ is given by a polycyclic presentation
described in Section \ref{se:poly_pres}.
First we show that the group $G_F$ can be
presented as a semidirect product of an abelian matrix group and $\MZ^n$.
Then we present the attack on the obtained presentation.

\subsection{Deduced semidirect product for $G_F$}
\label{se:deduced_pres}

Given a polycyclic presentation for $G_F$ constructed in Section \ref{se:poly_pres}
it is straightforward to find the numbers $m$ and $n$.
For the relations:
$$
g_{m+j}^{g_i} = g_{m+1}^{a_{ij1}}\ldots g_{m+n}^{a_{ijn}}
$$
we can define matrices $C_1, \ldots, C_m$:
$$
C_{i} = (a_{ijk})_{j=1,\ldots,n}^{k=1,\ldots,n}.
$$
Next we form a semidirect product $G$ of $\gp{C_1, \ldots, C_m}$ and $\MZ^n$
which is a set of pairs:
$$
\{(C, \ovs) \mid C \in \langle C_1, \ldots, C_m\rangle, \ovs \in \mathbb{Z}^n\}
$$
equipped with the multiplication given by
$$
(C, \ovs)\cdot(D, \ovt) = (CD, \ovs D + \ovt).
$$
Let $\{\ove_1,\ldots,\ove_n\}$ be the standard basis for $\MZ^n$.
It is easy to check the map $\tau\colon\{g_1,\ldots,g_{m+n}\} \to G$ given by:
$$
\tau(g_i) =
\begin{cases}
(C_i,\ovo) & \mbox{if } i\le m,\\
(E,\ove_j) & \mbox{if } i=m+j,\ 1 \le j\le n,\\
\end{cases}
$$
defines an isomorphism between $G_F$ and the constructed group.
Furthermore, given an element  $g = g_1^{e_1} \ldots g_n^{e_n}$
it requires polynomial time to find its $\tau$-image.

We also claim that given a pair $(C,\ovv)$ it requires polynomial time
to find a word $g$ such that $\tau(g)=(C,\ovv)$.
To convert $(C,\ovv)$ into a word in the generators $g_1, \ldots, g_{m+n}$
one can express $(C,\ovv)$ as a product:
$$
(C, \ovv) = (C_1, \ovo)^{a_{1}}\ldots (C_m, \ovo)^{a_{m}}
(E, \ove_1)^{a_{m+1}}\ldots (E, \ove_n)^{a_{m+n}},
$$
for some $a_1, \ldots, a_{n+m} \in \mathbb{Z}$,
in which case $g = g_1^{a_1}\ldots g_n^{a_n}g_{m+1}^{a_{m+1}}\ldots g_{m+n}^{a_{m+n}}$.
Clearly $(C, \ovv) = (C, \ovo)(E, \ovv)$. Therefore we have to solve two tasks.
First, we need to find $a_{1}, \ldots, a_{m}$ such that
$C = C_1^{a_{1}}\ldots C_m^{a_{m}}$
which can be done in polynomial time~\cite{babai1996multiplicative}.
Second, we need to find $a_{m+1}, \ldots, a_{m+n}$ such that
$\ovv = a_{m+1}\ove_1 + \ldots + a_{m+n}\ove_n$ which is obvious.

It follows from the discussion above that computational problems
for $G_F$ given by polycyclic presentation and by the deduced semidirect product
are polynomial time equivalent. Another important property of the computed
presentation is that the ring:
$$
K = \mathbb{Q}[C_1, \ldots, C_m]
$$
generated by
matrices $C_1, \ldots C_m$ is actually a field isomorphic to a subfield of $F$
(because $C_i$'s define the same action as $U_i$'s, but in a basis
$O_1,\ldots,O_n$).

\subsection{The attack}

In the deduced presentation of $G_F$
the system of conjugacy equations (\ref{eq:conds})
is equivalent to the following system
of equations with unknown $C\in K^\ast$ and $\ovv\in\MZ^n$:
\begin{equation}\label{eq:sys1}
\left\{
\begin{array}{lcl}
\ovv(E - B_1) + \ovt_1C & = & \ovt'_1, \\
& \vdots & \\
\ovv(E - B_{N_2}) + \ovt_{N_2}C & = &\ovt'_{N_2},
\end{array}
\right.
\end{equation}
where $(C, \ovv)$ represents $X$, $(B, \ovt_i)$ represents $b_i$,
$(B, \ovt'_i)$ represents $b'_i$ for $i = 1, \ldots, N_2$.

To solve the system (\ref{eq:sys1}) we
compute a basis $H_1, \ldots, H_l$ of the field $K$
as a vector space over $\mathbb{Q}$. Hence,
$$C = c_1H_1 + \ldots + c_lH_l$$
for some $c_1,\ldots,c_n\in\MQ$ and (\ref{eq:sys1}) can be rewritten as:
\begin{equation*}\label{eq:sys2}
\left\{
\begin{array}{lcl}
\ovv(E - B_1) + c_1\ovt_1H_1 + \ldots + c_l\ovt_1H_l & = & \ovt'_1, \\
& \vdots & \\
\ovv(E - B_{N_2}) + c_1\ovt_{N_2}H_1 + \ldots + c_l\ovt_{N_2}H_l & = &\ovt'_{N_2},
\end{array}
\right.
\end{equation*}
which is a system of linear equations over field $\mathbb{Q}$ with unknown
$\ovv=(v_1, \ldots, v_n)$ and $c_1, \ldots, c_l\in \MQ$.
The solution of this system provides us with the key $A'$.

We call this procedure as \textbf{FBA2}. The attack also was implemented in GAP
and tested on the same machine.
The table bellow contains results of our tests.

\begin{center}
\begin{tabular}{|c|c|c|c|c|c|}
\hline
\raisebox{-4ex}[0cm][0cm]{Polynomial} & \raisebox{-4ex}[0cm][0cm]{$h(G)$} & \multicolumn{2}{|c|}{\raisebox{-1.5ex}[0cm][0cm]{FBA2, $L=5$}}  & \multicolumn{2}{|c|}{\raisebox{-1.5ex}[0cm][0cm]{FBA2, $L=100$}} \\
 & & \multicolumn{2}{|c|}{}  & \multicolumn{2}{|c|}{} \\
\cline{3-6}
& & \raisebox{-1.5ex}[0cm][0cm]{Time} & Success & \raisebox{-1.5ex}[0cm][0cm]{Time} & Success\\
&  &    &  rate & & rate\\
\hline
 $x^2-x-1$      &  3 &   4.3 s & 100\% &   3.9 s & 100\% \\
 $x^5-x^3-1$    &  7 &   4.9 s & 100\% &   6.8 s & 100\% \\
 $x^7-x^3-1$    & 10 &   8.1 s & 100\% &  10.1 s & 100\% \\
 $x^9-7x^3-1$   & 14 &  34.0 s & 100\% &  47.7 s & 100\% \\
 $x^{11}-x^3-1$ & 16 &  20.9 s & 100\% &  26.4 s & 100\% \\
 $x^{15}-x-2$   & 22 & 528.2 s & 100\% & 761.3 s & 100\% \\
 $x^{20}-x-1$   & 30 & 164.6 s & 100\% & 208.2 s & 100\% \\
\hline
\end{tabular}
\end{center}

\section{Conclusion}

Our arguments show the following.
\begin{itemize}
\item
The groups of the form $U_F\ltimes \CO_F$ can not be used as platform groups
in the commutator key-establishment protocol.
\item
It is difficult to devise a successful length-based-attack and low
success rate does not mean much in terms of security.
\end{itemize}
Finally we want to point out that our attack does not eliminate all polycyclic groups
from consideration.

\end{document}